\documentclass[12pt]{amsart}

\usepackage{amsmath,amsthm}
\usepackage{amsfonts,amssymb}
\usepackage{accents}
\usepackage{enumerate}
\usepackage{accents,color}
\usepackage{graphicx}

\usepackage{mathrsfs}
\usepackage{mathtools 
}

\usepackage[
pdftex,hyperfootnotes]{hyperref}

\usepackage{nicematrix}

\usepackage[ansinew]{inputenc}
\usepackage[T1]{fontenc}

\usepackage{textcomp}
\usepackage{charter}
\usepackage[]{bera}
\usepackage[charter]{mathdesign}

\usepackage{longtable}
\usepackage{enumerate}
\allowdisplaybreaks

\newcommand{\lvt}{\left|\kern-1.35pt\left|\kern-1.3pt\left|}
\newcommand{\rvt}{\right|\kern-1.3pt\right|\kern-1.35pt\right|}

\newtheorem{thm}{Theorem}[section]
\newtheorem{cor}[thm]{Corollary}

\newtheorem{exam}[thm]{Example}

\theoremstyle{remark}

 \def\d{\mathrm{d}}

 \def\dd{{\operatorname d}} 

\def\OP{\pmb{\operatorname P}}
\def\OT{\pmb{\operatorname T}}

\def\OH{\pmb{\operatorname H}}
\def\OI{\pmb{\operatorname I}}

 \def\BB{{\mathbb B}}
 \def\CC{{\mathbb C}}
 
 \def\NN{{\mathbb N}}

 \def\RR{{\mathbb R}}

 \def\ZZ{{\mathbb Z}}

\def\lla{\langle{\kern-2.5pt}\langle} 
\def\rra{\rangle{\kern-2.5pt}\rangle}

\graphicspath{{./}}

\begin{document}
 
\title{
Markov Chains and Multiple Orthogonality}

\author[Branquinho]{Amílcar Branquinho}
\address{CMUC, Department of Mathematics, University of Coimbra, 3000-143 Coimbra, Portugal}
\email{ajplb@mat.uc.pt}

\author[Díaz]{Juan E. F. Díaz}
\address{CIDMA, Departamento de Matemática, Universidade de Aveiro, 3810-193 Aveiro, Portugal}
\email{juan.enri@ua.pt}

\author[Foulquié-Moreno]{Ana Foulquié-Moreno}
\address{CIDMA, Departamento de Matemática, Universidade de Aveiro, 3810-193 Aveiro, Portugal}
\email{foulquie@ua.pt}

\author[Mañas]{Manuel Mañas}
\address{Departamento de Física Teórica, Universidad Complutense de Madrid, Plaza Ciencias 1, 28040-Madrid, Spain \& Instituto de Ciencias Matematicas (ICMAT), Campus de Cantoblanco UAM, 28049-Madrid, Spain}
\email{manuel.manas@ucm.es}


\date{\today} 
\subjclass[2020]{42C05,33C45,33C47,60J10,60Gxx}

\keywords{Multiple orthogonal polynomials, Markov chains, stochastic matrices, Karlin--McGregor representation formula, 
Poincaré's theorem
.}
 
\begin{abstract}
In this work we survey on connections of Markov chains and the theory of multiple orthogonality.
Here we mainly concentrate on give a procedure to generate stochastic tetra diagonal Hessenberg matrices, coming from some specific families of multiple orthogonal, such as the ones of Jacobi--Piñeiro and Hypergeometric Lima--Loureiro.
We show that associated with a positive tetra diagonal nonnegative bounded Hessenberg matrix we can construct two stochastic tetra diagonal ones.
These two stochastic tridiagonal nonnegative Hessenberg matrices are shown to be, enlightened by the Poincaré theorem, limit transpose of each other.
\end{abstract} 
 
\maketitle

%

\section{Introduction}\label{sec:0}

Since the works of Karlin and McGregor (cf. \cite{KmcG,Karlin-McGregor}) it is well known that the theory of orthogonal polynomials plays a key role in understanding Random Walks and their continuous Birth and Death process.

These comes from the fact that the transition matrix of these Markov chains or processes admits a tridiagonal Jacobi representation (cf.~\cite{Gallager}), and by the Stone--Shohat--Favard theorem (cf. \cite{Simon}) is known that there is a sequence of orthogonal polynomials connected with them.
In fact, 
we can determine the so called spectral measure as an uniform limit of quotient of associated polynomials and the orthogonal ones (which is the famous Markov--Stieltjes convergence theorem).
These two sequences of orthogonal polynomials (that exists by the Favard Theorem) are defined by the three term recurrence relations governed by the stochastic transition matrix.

In this work we will explain how to apply this technicalities coming from the constructive approximation theory to describe the transition matrices that govern the discrete Markov chains.

In a series of works we try to understand Markov Chains that goes beyond simple Random Walks admitting connections far from near neighbors (cf.~\cite{nuestro2,espectral,nuestro1}).
In this case the transitions matrix is no more of Jacobi type and the discrete Markov chain gives rises to a Hessenberg structure. Here the model of orthogonal polynomials do not apply, nevertheless some biorthogonality can be used to modelize and give some insight on the procedure.
These biorthogonality comes from the ideas of the multiorthogonality presented in several monographs (cf.~\cite{Ismail,nikishin_sorokin,VanAssche2}).

In this work we will focus on the direct problem related to a discrete Markov chain that come from the generalization of simple Random Walks and is associated with a tetra diagonal transition stochastic matrix.
We establish a Karlin--McGregor representation theorem, we characterize the transient and recurrent situations, and present some examples (related with the so called Jacobi--Piñeiro and Hypergeometric multiple orthogonal polynomials) of this model that come from the theory of multiple orthogonal polynomials.
We notice that there are interesting applications of these to urn models can be find in~\cite{Grunbaum_Iglesia}.

The work is organized as:
In Section~\ref{sec:1} we give a quick overview of the important results of the theory of orthogonal polynomials, fundamental to understand the Random Walks model.
In Section~\ref{sec:2} we explain how the theory of multiple orthogonal polynomials apply to Hessenberg, bounded, positive, transitions matrices, to give stochastic matrices.
We concentrate on the tetra diagonal situations and find out two configurations for the transitions stochastic matrices.
In Section~\ref{sec:3} we present the Karlin--McGregor representation theorem and
apply a Poincaré theorem get a limit relations between the two transitions stochastic matrices previously obtained.

\section{Orthogonal polynomials}\label{sec:1}

In their paper on Random Walks (cf. \cite{Karlin-McGregor}), Karlin and McGregor consider a Jacobi, nonnegative, stochastic semi-infinite transition matrix of type
\begin{align}
\label{eq:stotransmatrix}
\OP & = 
\resizebox{.3\textwidth}{!}{$\displaystyle
\left[
\begin{NiceMatrix}[columns-width = 0.5cm]
q_0 & r_0& 0 & \Cdots{} & & \\
p_1 & q_1 & r_1 & \Ddots & & \\
0 & p_2& q_2& r_2 & & \\
\Vdots {} 
&\Ddots {} 
& \Ddots {} 
& \Ddots {}
& \Ddots {} 
& \Ddots {} 
\end{NiceMatrix}
\right]$}
\end{align}
with $(q_n) \subset [0,+\infty[$, and $(p_n)$, $(r_n) \subset ]0,+\infty[ $, where
\begin{align} \label{eq:cond_estoc}
q_0 + r_0 = 1 , && 
p_n + q_n + r_n = 1, && n \in \NN ,
\end{align}
As an example we have the Chebychev transition matrix given by
\begin{align*}
\OT
 =
\resizebox{.25\textwidth}{!}{$\displaystyle
\left[\begin{NiceMatrix}
0 & 1 & 0 & \Cdots & \\
\frac1 2 & 0 & \frac 12 & \Ddots & \\
 0 & \frac 1 2 & 0 & \frac 1 2 & \\
\Vdots & \Ddots & \Ddots & \Ddots & \Ddots
 \end{NiceMatrix} \right]$} .
\end{align*}
In fact, the Chebychev polynomials are defined as
\begin{align*}
T_n (x) = \cos (n \vartheta) && \text{with} && x = \cos (\vartheta) , && n \in \NN \cup \{ 0 \} ,
\end{align*}
or equivalently by, $T_{-1} (x) = 0$, 
$T_0 (x) = 1$, 
and
\begin{align*}
x \, T_n (x) & = \frac 12 T_{n+1} (x) + \frac 12 T_{n-1} (x) , && n \in \NN .
\end{align*}
We can condense all this information in the matrix equation
\begin{align}
\label{eq:operador_chebychev}
 x \, \mathcal T = \OT \, \mathcal T && \text{where} && \mathcal T = \left[ \begin{NiceMatrix} T_0 (x) & T_1 (x) & \Cdots & T_n (x) & \Cdots \end{NiceMatrix} \right]^\top .
\end{align}
But we know that,
\begin{align*}
\int_{-1}^1 T_n (x) \, T_m (x) \, \frac{1}{\pi \, \sqrt{1-x^2}} \, \dd x = \frac 1 2 \, \delta_{n,m} , && n,m \in \NN ,
\end{align*}
which asserts that $\big\{ T_n \big\}$ is a sequence of polynomials, orthogonal with respect to the weight function, $w (x) = \frac 1 {\pi \, \sqrt{1-x^2}}$, defined on $[-1,1]$.
The orthogonality property enables us to describe the \emph{evolution} of the transition matrix. In fact, from~\eqref{eq:operador_chebychev} we get
\begin{align*}
x^k \, \mathcal T = \OT^k \, \mathcal T , && k \in \NN ,
\end{align*}
and so, $\OT^k = \left[ \begin{NiceMatrix} \OT_{n,m}^k \end{NiceMatrix} \right]$, where
\begin{align*}
\OT_{n,m}^k =
\frac{\displaystyle \int_{-1}^1 x^k \, T_n (x) \, T_m (x) \frac{1}{\pi \, \sqrt{1-x^2}} \, \dd x }
{\displaystyle \int_{-1}^1 x^k \, T_n^2 (x) \frac{1}{\pi \, \sqrt{1-x^2}} \, \dd x } , 
 && n,m, k \in \NN \cup \{ 0 \} .
\end{align*}
Now, we will see that the Chebychev orthogonal polynomial sequences 
is a model to describe Jacobi transition matrices.
In fact, associated with the transition matrix, $\OP$, in~\eqref{eq:stotransmatrix} we can define a sequence of polynomials $\big\{ P_n \big\}$ by means of
\begin{align}
\label{eq:operacional}
 x \, \mathcal P = \OP \, \mathcal P && \text{where} && \mathcal P
 = \left[ \begin{NiceMatrix}
 P_0 (x) & P_1 (x) & \Cdots & P_n (x) & \Cdots \end{NiceMatrix} \right]^\top ,
\end{align}
and so, taking, $P_{-1} (x) = 0$, $P_0 (x) = 1$
\begin{align}
\label{eq:ttrr}
x \, P_n (x) & = r_n \, P_{n+1} (x) + q_n \, P_{n} (x) + p_n \, P_{n-1} (x) , && n \in \NN .
\end{align}
We now state one of the main theorems in the theory of orthogonal polynomials.

\begin{thm}[Stone--Shohat--Favard]
Let $\big\{ P_n \big\}$ be a sequence of polynomials defined by~\eqref{eq:ttrr} with initial conditions $P_{-1} (x) = 0$, $P_0 (x) = 1$.
If $p_n,q_n,r_n$ satisfies~\eqref{eq:cond_estoc} then, there exists a 
positive measure, $\dd \mu$, defined on $[-1,1]$ such that $\big\{ P_n \big\}$ is a sequence of polynomials with respect to $\dd \mu$, i.e.
\begin{align*}
\int_{-1}^1 P_n (x) \, P_m (x) \, \dd \mu (x) = \kappa_n \, \delta_{n,m} , && n,m \in \NN \cup \{ 0 \} ,
\end{align*}
with $\kappa_n >0 $, $n \in \NN \cup \{ 0 \} $.
\end{thm}

The converse of this theorem is an easy consequence of the Fourier expansion for $x \, P_n$ in terms of the elements of the sequence $\big\{ P_n \big\}$.

Although this theorem gives the existence of a
measure describing the stochastic transition matrix, do not explain how to determine it from the data. Nevertheless, we can determine the sequence of moments for that measure from the~\eqref{eq:ttrr} (cf. for instance \cite{Chihara}).

In the works, \cite{KmcG,Karlin-McGregor}, the authors explain that the moment problem associated with this type of transition matrices is determined, and can be represented by its absolutely continuous parts up to a discrete part (complete determined by the data).
This is a key fact, as from the Markov--Stieltjes theorem (in its weak form) we get the spectral measure of~$\OP$.

\begin{thm}[Markov--Stieltjes]
Let $\big\{ P_n \big\}$ be a sequence of polynomials with respect to 
positive measure defined on $[-1,1]$; then
\begin{align*}
\lim_{n \to \infty} \frac{P_{n-1}^{(1)} (z)}{P_n (z)} = \int_{-1}^1 \frac{\dd \mu (x)}{z-x} , 
\end{align*}
holds uniformly on compact sets of $\CC \setminus [-1,1]$.
\end{thm}

The polynomial sequence, $\big\{ P^{(1)}_n \big\}$, defined by
\begin{align*}
P_{n-1}^{(1)} (z) = 
 \int_{-1}^1 \frac{P_{n+1} (z) - P_{n+1} (x)}{z-x} \, \dd \mu (x) , && n \in \mathbb N ,
\end{align*}
supposing that $\displaystyle \int_{-1}^1 \dd \mu (x) = 1$, 
is the associated polynomials sequence of~$\big\{ P_n \big\}$ and~$\mu$ on $[-1,1]$.
It is important to note that the Markov--Stieltjes function,
\begin{align}
\label{eq:Markov_Stieltjes}
 S (z) = \int_{-1}^1 \frac{\dd \mu (x)} 
 {z-x} 
 , && \text{analytic in} && \CC \setminus [-1,1] ,
\end{align}
is the spectral measure of the operator defined by the stochastic transition matrix $\OP$.
Moreover, we can recover the absolutely continuous part of the measure $\mu$ by using the Sokhotski--Plemelj formula (cf.~\cite{Deift}), as
\begin{align*}
\big( S(z) \big)_+ - \big( S(z) \big)_- = 2 \pi \, \pmb{\operatorname i} \, w (x) , && \text{where} &&
\big( S(z) \big)_\pm = \lim_{\varepsilon \to 0^\pm} f (x + \pmb{\operatorname i} \, \varepsilon) .
\end{align*} 
In fact, this function can be viewed as a complex measure of orthogonality for $\big\{ P_n \big\}$ on a disc $D \subset \CC \setminus [-1,1]$ with positive orientation:
\begin{align*}
 & \frac{1}{2 \pi \, \pmb{\operatorname i}} \, \int_{\partial D} P_{n}(z) \, P_{m}(z) \, S (z) \, \operatorname d z \\
 & \phantom{olaolaola} = \frac{1}{2 \pi \, \pmb{\operatorname i}} \int_{\partial D} P_{n}(z) \, P_{m}(z) \int_{-1}^1 \frac{\dd \mu (x)
 }{z-x} \operatorname d z && \text{by~\eqref{eq:Markov_Stieltjes}} \\
 & \phantom{olaolaola} = \int_{-1}^1 
 \dd \mu (x) \frac{1}{2\pi \, \pmb{\operatorname i}} \int_{\partial D} \frac{P_{n}(z) \, P_{m}(z)}{z-x} \operatorname d z && \text{(Fubini Theorem)}
 \\ 
& \phantom{olaolaola} = \int_{-1}^1 P_{n}(x) \, P_{m}(x) \, \dd \mu (x) 
. && \text{(Cauchy Theorem)} 
\end{align*}
We can say that the spectral measure of a Jacobi operators connects the theories of orthogonal polynomials with the one of operators.

In the proof of Markov--Stieltjes theorem heavily use that the zeros of orthogonal polynomials $P_n$ are real, simple, lies on $[-1,1]$ and have the interlacing property.
Note that the zeros of the orthogonal polynomials $P_n$ are for each $n \in \NN$ the eigenvalues of the principal truncated matrix of~$\OP$.
All of these follows from the fact that for stochastic positive matrices,~$\OP$, there exists a real number $c$ such that $\OP - c \, \OI$ is oscillatory (corollary of the Gantmacher and Krein criteria for a Jacobi matrix, cf.~\cite[Theorem 1.7]{espectral}). For a nice account on oscillatory matrices cf. for instance~\cite{Fallat-Johnson,Gantmacher-Krein}.

\begin{thm}
The spectral measure of the operator $\OP$ in~\eqref{eq:stotransmatrix} admits the representation
\begin{align*}
S (z) = \ell_0^\top \, \big( z \, \OI - \OP \big)^{-1} \, \ell_0 = \sum_{n = 0}^\infty \frac{\ell_0^\top \, \OP^n \, \ell_0}{z^{n+1}}, && |z| > 1 ,
\end{align*}
where
$ \ell_{0}^{\top} = \left[ \begin{NiceMatrix} 1 & 0 & \Cdots \end{NiceMatrix} \right]$.
\end{thm}

This result is a easy consequence of the representation of the moments
$\displaystyle
w_{n} 
 = \int z^{n} \, S (z) \, \operatorname d z $, $n \in \NN$, the three term recurrence relation for $\big\{ P_n \big\}$~\eqref{eq:operacional} and the representation of the spectral measure~\eqref{eq:Markov_Stieltjes}.

Now, we can describe the \emph{evolution} of the transition matrix
$\OP^k = \left[ \begin{NiceMatrix} \OP_{n,m}^k \end{NiceMatrix} \right]$,~by
\begin{align*}
\OP_{n,m}^k =
\frac{\displaystyle \int_{-1}^1 x^k \, P_n (x) \, P_m (x) \, \dd \mu (x) }
{\displaystyle \int_{-1}^1 x^k \, P_n^2 (x) \, \dd \mu (x) } , 
 && n,m, k \in \NN .
\end{align*}

\section{Multiple orthogonal polynomials}\label{sec:2}

In this section we show how the multiple orthogonal polynomials theory can be applied to construct stochastic matrices. Here we follow mainly the works \cite{afm,nuestro1}.

When we deal with monic multiple orthogonal on the stepline we have to consider bounded, non-negative, Hessenberg matrices of type
\begin{align*}
\OH
 =
 \resizebox{.3\textwidth}{!}{$\displaystyle
 \left[ \begin{NiceMatrix}
 c_0 & 1 & 0 & \Cdots & & \\[-5pt]
 b_1 & c_1 & 1 & \Ddots& & \\ 
 a_2 & b_2 & c_2 & 1 &&\\ 
 0& a_3 & b_3 & c_3 & 1 & \\
 \Vdots&\Ddots&\Ddots&\Ddots&\Ddots&\Ddots
 \end{NiceMatrix} \right]$}
 && \text{with} && a_n > 0 , && n \in \NN 
\end{align*}
(we restrict ourself with no loss of generality to the tetra diagonal case).
Associated with~$\OH $ we formally define three sequences of polynomials, $\big\{ B_n \big\}$,
$\big\{ A_n^1 \big\}$ and $\big\{ A_n^2 \big\}$ by means of
\begin{align}
\label{eq:tetrarrB}
 x \, \mathcal B & = \OH \, \mathcal B && \text{where} && \mathcal B = \left[ \begin{NiceMatrix} B_0 (x) & B_1 (x) & \Cdots & B_n (x) & \Cdots \end{NiceMatrix} \right]^\top , \\
\label{eq:tetrarrA1} 
 x \, \mathcal A^1 & = \mathcal A^1 \, \OH && \text{where} && \mathcal A^1 = \left[ \begin{NiceMatrix} A_0^1 (x) & A_1^1 (x) & \Cdots & A_n^1 (x) & \Cdots \end{NiceMatrix} \right] , \\
\label{eq:tetrarrA2} 
\hspace{-.5cm} x \, \mathcal A^2 & = \mathcal A^2 \, \OH && \text{where} && \mathcal A^2 = \left[ \begin{NiceMatrix} A_0^2 (x) & A_1^2 (x) & \Cdots & A_n^2 (x) & \Cdots \end{NiceMatrix} \right] , 
\end{align}
with initial conditions
\begin{align*}
B_0 (x) = 1 , && A_{-1}^1 (x) = 0 , && A_0^1 (x) = 1 , && A_{-1}^2 (x) = 1 , && A_0^2 (x) = 0 .
\end{align*}
Let us consider a couple of weights $(w_1,w_2)$, the semi-infinite vectors of monomials
\begin{align*}
\chi & = \left[ \begin{NiceMatrix}
 1 & x & x^2 & \Cdots 
 \end{NiceMatrix} \right]^\top, \\
\chi_1 & =\left[ \begin{NiceMatrix}
 1 & 0 & x & 0 & x^2 & 
 \Cdots 
 \end{NiceMatrix} \right]^\top, 
 \\
\chi_2 & = \left[ \begin{NiceMatrix}
 0 & 1 & 0 & x & 0 & x^2 &
 \Cdots
 \end{NiceMatrix} \right]^\top,
\end{align*}
and the following vector of undressed linear forms
\begin{align*}
 \xi &:= \chi_1 w_1 + \chi_2 w_2
 =\left[ \begin{NiceMatrix}
 w_1 & w_2 & xw_1 & xw_2 & x^2w_1 & x^2w_2 &
 \Cdots
 \end{NiceMatrix} \right]^\top. 
\end{align*}
Given a Borel positive measure $\dd \mu$ with support on the closed interval $\Delta \subset \RR$, our moment matrix is 
 \begin{align*}
 g =\int_\Delta \chi(x) \, \big( \xi(x) \big)^\top \, \dd \mu(x).
 \end{align*}
 The Gauss--Borel factorization of the moment matrix $g$ (that exists and is unique whenever all its principal minors are non-singular) is the problem of finding the solution of
 \begin{align*}
 g &=S^{-1} \, H \, \tilde S^{-\top}, 
 \end{align*}
 with $S$, $\tilde S$ lower unitriangular semi-infinite matrices
 \begin{align*}
 S&=
 \resizebox{.3\textwidth}{!}{$\displaystyle
 \left[ \begin{NiceMatrix}[columns-width = auto]
 1 & 0 & \Cdots & \\
 S_{1,0 } & 1& \Ddots & \\
 S_{2,0} & S_{2,1} & \Ddots & \\
 \Vdots & \Ddots & \Ddots &
 \end{NiceMatrix}\right]$}, & 
 \tilde S&=
 \resizebox{.3\textwidth}{!}{$\displaystyle
 \left[ \begin{NiceMatrix}[columns-width = auto]
 1 & 0 &\Cdots &\\
 \tilde S_{1,0 } & 1&\Ddots&\\
 \tilde S_{2,0} & \tilde S_{2,1} & \Ddots &\\
 \Vdots & \Ddots& \Ddots& 
 \end{NiceMatrix} \right]$} , 
 \end{align*}
 and $H$ a invertible, semi-infinite diagonal matrix 
$ H=\operatorname{diag}
 \left[ \begin{NiceMatrix}
 H_0 & H_1 & \Cdots 
 \end{NiceMatrix} {} \right]$.
The type~II multiple orthogonal polynomials, $\big\{ B_n \big\}$, and of type~I linear forms, $\big\{ Q_n \big\}$, with
 \begin{align*}
 Q_n (x) = A^1_n (x) \, w_1 (x) + A^2_n (x) \, w_2 (x) , && n \in \NN \cup \{ 0 \} ,
 \end{align*}
are defined by
 \begin{align*}
 \mathcal B & = S \, \chi, &&
 \mathcal Q = \left[ \begin{NiceMatrix} Q_0 (x) & Q_1 (x) & \Cdots & Q_n (x) & \Cdots \end{NiceMatrix} \right]
= H^{-1} \, \tilde S \, \xi , 
 \end{align*}
They fulfill the following type I
 \begin{align*}
 \int_{\Delta} x^{j} \, Q_n (x)
 \, \dd \mu (x) =0,
 && j\in\{0,\ldots, n - 1\} ,
 \end{align*}
 and type~II orthogonality relations
\begin{align*}
 \int_{\Delta} B_{n}(x) \, w_{a} (x) \, x^{j_a} \, \dd \mu (x) =0, && j_a=0,\ldots, \Big\lfloor \frac{n - a } 2 \Big\rfloor , && a=1,2 .
\end{align*}
We also have the following multiple biorthogonality relations 
 \begin{align*}
 \int_\Delta B_{m}(x) \, Q_{k}(x) \, \dd \mu(x) = \delta_{m,k}, && m,k \in \NN \cup \{ 0 \} .
 \end{align*}
We can also show that the matrix $\OH$ associated with these families $\big\{ B_n \big\}$ and~$\big\{ Q_m \big\}$ is defined by
\begin{align*}
\OH := S \, \Lambda \, S^{-1} 
 && \text{where} &&
 \Lambda 
 =
 \resizebox{.25\textwidth}{!}{$\displaystyle
 \left[ \begin{NiceMatrix}
 0 & 1 & 0 & \Cdots \\
 0 & 0 & 1 & \Ddots \\
 \Vdots & \Ddots & \Ddots & \Ddots
 \end{NiceMatrix}\right]$}.
\end{align*}
When we are in presence of an 
algebraic Chebychev system $(w_1 \, \dd \mu , w_2 \, \dd\mu)$ (AT-systems for short)
we can see that the zeros of the type II polynomials or type I linear forms are contained on $\BB (0, \|\OH \|) = \big\{ x \in \CC : | x | < \|\OH\| \big\}$ (where $\| \OH \|$ is the usual operator norm of $\OH$) and have the interlacing property (cf. \cite{Ismail,nikishin_sorokin,VanAssche2}).
These enables us to prove (cf. \cite[Lemma 2 and~3]{nuestro1})~that
\begin{align}
\label{eq:positividad}
Q_{n} (\lambda) > 0, && B_{n} (\lambda) > 0, && \lambda \geq \| \OH \| , && n \in \NN .
\end{align}
Hence, taking $x = \| \OH \| \, t$ in~\eqref{eq:tetrarrB} we get, 
\begin{multline*}
a_n \, B_{n-2} (\| \OH \| \, t) + b_{n} \, B_{n-1} (\| \OH \| \, t) + c_n \, B_n (\| \OH \| \, t) + B_{n+1} (\| \OH \| \, t) 
 \\ = \| \OH \| \, t \, B_n (\| \OH \| \, t)
\end{multline*}
and dividing by $\| \OH \|^{n+1}$ we arrive to
\begin{align*}
\tilde a_n \, \tilde B_{n-2} ( t)
 + \tilde b_n \, \tilde B_{n-1} ( t)
 + \tilde c_n \, \tilde B_n ( t)
 + \tilde B_{n+1} (t) = t \, \tilde B_n ( t) ,
\end{align*}
with
\begin{align*}
\tilde a_n = \frac{a_n}{\|\OH\|^3}, &&
\tilde b_n = \frac{b_n}{\|\OH\|^2} , &&
\tilde c_n = \frac{c_n}{\|\OH\|} , &&
\tilde B_n ( t) = \frac{B_{n} (\| \OH \| \, t )}{\| \OH \|^n} ,
\end{align*}
We can see from here that, when $\OH$ is bounded, we can take a similar operator $ \tilde \OH$ with norm one, i.e.
\begin{align*}
\tilde \OH = \sigma^{-1} \, \frac 1 \lambda \, \OH \, \sigma,
&& \text{with} &&
\sigma = \operatorname{diag}
\left[
\begin{NiceMatrix}
1 & \| \OH \| & \Cdots & \| \OH \|^n & \Cdots
\end{NiceMatrix}
\right] 
.
\end{align*}
Note that the new eigenvectors are
\begin{align*}
 \tilde{\mathcal B} & = \left[ \begin{NiceMatrix} \tilde B_0 (x) & \tilde B_1 (x) & \Cdots & \tilde B_n (x) & \Cdots \end{NiceMatrix} \right]^\top , && \text{(right eigenvector)} \\
 \tilde{\mathcal Q} & = \left[ \begin{NiceMatrix} \tilde Q_0 (x) & \tilde Q_1 (x) & \Cdots & \tilde Q_n (x) & \Cdots \end{NiceMatrix} \right] , && \text{(left eigenvector)}
\end{align*} 
with $\displaystyle \tilde Q_n (x) = \frac{Q_{n} (\| \OH \| \, t )}{\| \OH \|^n} $.

From now on we consider, with no loss of generality, norm one Hessenberg matrix, $\OH$. 

Now we will see how to construct by similarity a stochastic transition matrix associated with $\OH$.
The procedure is constructive.

As we know from~\eqref{eq:positividad} we can multiply equation~\eqref{eq:tetrarrB} by $B_n (1)$, (respectively, the linear combination of~\eqref{eq:tetrarrA1} and~\eqref{eq:tetrarrA2} multiplied by $w_1$ and~$w_2$ by $Q_n (1)$), to get
\begin{align*}
\hat a_n \, \hat B_{n-2} (x)
 + \hat b_n \, \hat B_{n-1} (x)
 + c_n \, \hat B_n (x)
 + \hat d_n \, \hat B_{n+1} (x) = x \, \hat B_n (x) , \\
\alpha_n \, \check Q_{n-1} (x)
 + c_n \, \check Q_{n} (x)
 + \beta_n \, \check Q_{n+1} (x)
 + \gamma_n \, \check Q_{n+2} (x) = x \, \check Q_n (x) , 
\end{align*}
with
\begin{align*}
\hat a_n & = \frac{a_n \, B_{n-2} (1)}{B_{n} (1)} , &
\hat b_n & = \frac{b_n \, B_{n-1} (1)}{B_{n} (1)} , &
\hat d_n & = \frac{B_{n+1} (1)}{B_{n} (1)} , &
\hat B_n ( t) & = \frac{B_{n} (x)}{B_{n} (1)} , \\
\alpha_n & = \frac{Q_{n-1} (1)}{B_{n} (1)} , &
\beta_n & = \frac{b_{n+1} \, Q_{n+1} (1)}{Q_{n} (1)} , &
\gamma_n & = \frac{a_{n+2} \, Q_{n+2} (1)}{Q_{n} (1)} , &
\check Q_n (x) & = \frac{Q_{n} (x)}{Q_{n} (1)} .
\end{align*}
By construction we have
\begin{align*}
\hat a_n + \hat b_n + c_n + \hat d_n = 1 , && \text{as well as,} &&
\alpha_n + c_n + \beta_n + \gamma_n = 1 .
\end{align*}
Hence, we are in presence of two stochastic matrices associated with $\OH$
\begin{align}
\label{eq:estocastica2}
\hat \OH & = \sigma_{II}^{-1} \, \OH \, \sigma_{II} ,
&& \text{with} &&
\sigma_{II} = \operatorname{diag}
\left[
\begin{NiceMatrix}
1 & B_1 (1) & \Cdots & B_n (1) & \Cdots
\end{NiceMatrix}
\right] , \\
\label{eq:estocastica1}
\check \OH & = \sigma_{I}^{-1} \, \OH^\top \, \sigma_{I} ,
&& \text{with} &&
\sigma_{I} = \operatorname{diag}
\left[
\begin{NiceMatrix}
1 & Q_1 (1) & \Cdots & Q_n (1) & \Cdots
\end{NiceMatrix}
\right] .
\end{align}
The right eigenvectors for these matrices are
\begin{align*}
 \hat{\mathcal B} & = \left[ \begin{NiceMatrix} \hat B_0 (x) & \hat B_1 (x) & \Cdots & \hat B_n (x) & \Cdots \end{NiceMatrix} \right]^\top , \\
 \check{\mathcal Q} & = \left[ \begin{NiceMatrix} \check Q_0 (x) & \check Q_1 (x) & \Cdots & \check Q_n (x) & \Cdots \end{NiceMatrix} \right]^\top ,
\end{align*}
that is, 
$\displaystyle
x \, \hat{\mathcal B} = \hat \OH \, \hat{\mathcal B} $,
as well as, 
$\displaystyle
x \, \check{\mathcal Q} = \check \OH \, \check{\mathcal Q} $.

\begin{exam}[Jacobi--Piñeiro]\label{example:JP}
Consider the system of measures on $[0, 1]$:
\begin{align*}
w_1 (x) = x^{\alpha_1} (1-x)^{\alpha_0} \, \dd x , && w_2 (x) = x^{\alpha_2} (1-x)^{\alpha_0} \, \dd x,
\end{align*}
where
$\alpha_0 , \alpha_1 , \alpha_2 > -1$, $\alpha_1 - \alpha_2 \not \in \ZZ$. 
\end{exam}

The type~II Jacobi--Piñeiro multiple orthogonal polynomials are defined by a Rodrigues type formula (cf. \cite{abv})
\begin{multline*}
B_{n_1+ n} (x) = \frac{(1-x)^{-\alpha_0}}{M_{n_1 + n }} 
\Big( x^{-\alpha_2} \dfrac{\dd^{n} }{\dd x^{n}} x^{n + \alpha_2} \Big)
\Big( x^{-\alpha_1} \dfrac{\dd^{n_1} }{\dd x^{n_1}} x^{n_1 + \alpha_1} \Big) (1-x)^{n_1+n + \alpha_0} ,
\end{multline*}
with
$\displaystyle
M_{n_1 + n} = (-1)^{n_1+n} \frac
{ \Gamma (2 n_1 + n +\alpha_1+\alpha_0 + 1) \, \Gamma (n_1 + 2 n +\alpha_2+\alpha_0 + 1) }
{ \Gamma (n_1 + n +\alpha_1+\alpha_0 + 1) \, \Gamma (n_1 + n +\alpha_2+\alpha_0 + 1) }$, 
 and $n_1 = n$  or  $n_1 = n+1$, $n \in \NN$. 
We can write
\begin{align*}
a_{n+2} & = \lambda_{3n+2} \lambda_{3n+4} \lambda_{3n+6} \\
b_{n+1} & = \lambda_{3n+1} \lambda_{3n+3} + \lambda_{3n+2} \lambda_{3n+3} + \lambda_{3n+2} \lambda_{3n+4} \\
c_n & = \lambda_{3n} + \lambda_{3n+1} + \lambda_{3n+2} ,
\end{align*}
where for all $n \in \NN \cup \{ 0 \}$, we have (cf.~\cite{Aptekarev_Kaliaguine_VanIseghem})
\begin{align*}
\lambda_{6n+1} & = \frac{(2n + 1 + \alpha_1 + \alpha_0) \,
(2n+1 + \alpha_2 + \alpha_0) \,
(n+1 + \alpha_1)}
{(3n+1+\alpha_1 + \alpha_0) \,
(3n + 2 + \alpha_1 + \alpha_0) \,
(3n + 1 + \alpha_2 +\alpha_0)} , \\
\lambda_{6n+2} & =
\frac{(2n +1+\alpha_2 +\alpha_0) \, 
(2n +1+\alpha_0) \,
(n +\alpha_2 -\alpha_1)}
{(3n +1+\alpha_2 +\alpha_0) \,
(3n +2+\alpha_2 +\alpha_0) \,
(3n +2+\alpha_1 +\alpha_0)} , \\
\lambda_{6n+3} & =
\frac{(2n +2+1 +\alpha_0) \,
(2n +1+\alpha_0) \,
(n +1+\alpha_1 -\alpha_2)}
{(3n +2+\alpha_1 +\alpha_0) \,
(3n +3+\alpha_1 +\alpha_0) \,
(3n +2+\alpha_2 +\alpha_0)} , \\
\lambda_{6n+4} & =
\frac{(n +1+\alpha_2) \,
(2n +2+\alpha_2 +\alpha_0) \,
(2n +2+\alpha_1 +\alpha_0)}
{(3n +2+\alpha_2 +\alpha_0) \,
(3n +3+\alpha_2 +\alpha_0) \,
(3n +3+\alpha_1 +\alpha_0)} , \\
\lambda_{6n+5} & =
\frac{(n +1) \,
(2n +2+\alpha_2 +\alpha_0) \,
(2n +2+\alpha_0)}
{(3n +3+\alpha_2 +\alpha_0) \,
(3n +3+\alpha_1 +\alpha_0) \,
(3n +4+\alpha_1 +\alpha_0)} , \\
\lambda_{6n} & =
\frac{n \,
(2n+1+ \alpha_1 + \alpha_0) \,
(2n+ \alpha_0)}
{(3n + \alpha_2 + \alpha_0) \,
(3n +1+\alpha_2 +\alpha_0) \,
(3n +1+\alpha_1 +\alpha_0)} .
\end{align*}
In \cite{nuestro1} we have proved that the coefficients of the associated Hessenberg matrix are positive if and only if
the Jacobi--Piñeiro parameters, $\alpha_1 , \alpha_2$ are such that $|\alpha_1 - \alpha_2| < 1$.
But, the values of the type II, Jacobi--Piñeiro multiple orthogonal polynomials, at~$1$ are
\begin{align} \label{eq:Ben1}
B_{n_1 + n}(1)
 = \frac{(\gamma +1)_{n_1+n}}{(\alpha +\gamma +n_1+n+1)_{n_1} (\beta +\gamma +n_1+n+1)_n},
\end{align}
with  $n_1 = n$  or  $n_1 = n+1$,  $n \in \NN$. 
We also have
\begin{align}
\label{eq:limit}
\lim_{n\to \infty} {a_n} = \kappa^3 , && \lim_{n\to \infty} {b_n} = 3 \kappa^2 , && \lim_{n\to \infty} {c_n} = 3 \kappa , && \text{with} &&
\kappa = \frac 4{27}.
\end{align}
Now, we are able to construct, using~\eqref{eq:estocastica2}, the stochastic transition matrix,~$\hat \OH$, for the Jacobi--Piñeiro polynomials. To construct $\check \OH$, we need a close expression for the type I multiple orthogonal polynomials. In \cite{nuestro1} we got it for all $n \in \NN$,
\begin{align*}
 & \begin{multlined}[t][0.9\textwidth]
A_{2n}^1 (x) =
A_{2n}^2 (x)=
\frac{(2n+ \alpha + \gamma)_{n} \, (2n + \beta + \gamma)_{n}}{(n-1)! \, \Gamma( 2 n + \gamma)}
\\ \times \sum_{j= 0}^{n-1} (-1)^j \binom{n - 1}{j} \frac{\Gamma( 3 n - (j + 1) + \alpha + \gamma)}{\Gamma(n - j + \alpha) \, (\alpha - \beta - j)_n} \, x^{n - 1 - j} ,
\end{multlined}
\\
& \begin{multlined}[t][.9\textwidth]
A_{2n-1}^1 (x) =
\frac{\Gamma (3 n -1 + \alpha + \gamma) \, ( 2 n - 1+ \beta + \gamma)_n}{(n-1)! \, \Gamma ( 2 n - 1+\gamma )}
\\ \times \sum_{j=0}^{n-1} (-1)^{j} \binom{n-1}{ j}
\frac{(2 n - 1+\alpha + \gamma )_{n-1-j}}{\Gamma (\alpha + n - j) \, (\alpha - \beta + 1 - j)_{n-1}} \, x^{n -1 - j} ,
\end{multlined} \\
 & \begin{multlined}[t][.9\textwidth]
A_{2n+1}^2 (x) =
(-1)^{n-1}\frac{\Gamma (3 n + 1 + \beta + \gamma) \, (2 n + 1 + \alpha + \gamma)_{n+1}}{(n-1)! \, \Gamma(2 n + 1 + \gamma)}
\\ \times \sum_{j=0}^{n-1} (-1)^j \binom{n - 1}{ j}
\frac{( 2 n + 1+\beta +\gamma )_{n-1-j}}{\Gamma (\beta + n - j) \, (\alpha - \beta - (n - 1) + j)_{n+1}} \, x^{n - 1 - j} .
\end{multlined}
\end{align*}
For $\alpha_1 = -\frac{1}{4}$, $\alpha_2 = \alpha_0= -\frac{1}{2}$, 
we get, using~\eqref{eq:estocastica2} and~\eqref{eq:estocastica1}, the following approximate transition matrix in decimal form with a precision of four significant digits (now is not possible to find closed rational expressions and several sums involving the Euler Gamma function are required)
\begin{align*}
\check \OH
 \approx
 \resizebox{.44\textwidth}{!}{$\displaystyle
 \left[ \begin{NiceMatrix}[columns-width = 0.5cm]
0.6000& 0.2531& 0.1469& 0 & 0 &0&0& 
 \Cdots\\ 
0.4215&0.3167& 0.2419& 0.0199 & 0 &0&0& 
 \Ddots \\ 
0& 0.2760& 0.4657& 0.2036&0.0547&0&0& 
 \Ddots\\
0& 0 &0.3223&0.4176& 0.2341&0.0260&0& 
 \Ddots\\
0&0&0& 0.2826&0.4586&0.2110&0.0478&
 \Ddots\\
\Vdots&\Ddots&\Ddots&\Ddots&\Ddots&\Ddots&\Ddots&\Ddots&\Ddots
 \end{NiceMatrix}\right] $} , &&
\hat \OH
 \approx
 \resizebox{.385\textwidth}{!}{$ 
 \left[ \begin{NiceMatrix}
0.6000& 0.4000 & 0 & 0 & 0 &0&0& \Cdots \\
0.2667 & 0.3167 & 0.4167 & 0 & 0 &0&0& \Ddots \\
0.1026& 0.1603& 0.4657& 0.2715 &0&0&0& \Ddots \\
0& 0.0156 & 0.2417 &0.4176& 0.3250 &0&0& \Ddots \\
0&0&0.0565& 0.2035&0.4586&0.2814 &0&\Ddots \\
0&0&0&0.0250&0.2336& 0.4289 & 0.3125 & \Ddots \\
\Vdots&\Ddots&\Ddots&\Ddots&\Ddots&\Ddots&\Ddots&\Ddots
\end{NiceMatrix}\right]$} .
\end{align*}
As we can see these matrices are not transpose of each other.
We will try to give some light on this behavior in the next sections.

\section{Karlin--McGregor representation theorem}\label{sec:3}

We now extend the results of Karlin and McGregor concerning tridiagonal stochastic matrices~\cite{KmcG} to the multi-diagonal situation of multiple orthogonal polynomials.

We have seen in Section~\ref{sec:2}
that certain sets of multiple orthogonal polynomials give two families of stochastic matrices $\hat \OH$ and $\check \OH$. Here $\hat \OH$ models a 
Markov chain
with allowed jumps backward farther than near neighbors, and $\check \OH$ models a Random Walk with allowed jumps forward farther than near neighbors.

\begin{thm}[Karlin--McGregor representation formula]
\label{teo:KMcG}
Let $(w_1 \dd \mu,w_2 \dd \mu)$ be an AT-system.
Then, for a 
Markov chain with transition
matrix $\hat \OH $ (respectively, $\check \OH$), the 
probability, after~$r$ 
steps from state~$n$ to state $m$ is given~by
\begin{align*}
\hat \OH_{nm}^r & = \frac{B_{m}(1)}{B_{n}(1)} \int_\Delta x^r \, B_{n}(x) \, Q_{m}(x) \, \d \mu(x), 
 \\
\check \OH_{nm}^r & = \frac{ Q_{m}(1)}{Q_{n}(1)} \int_\Delta x^r \, B_{m}(x) \, Q_{n}(x) \, \d \mu(x) .
\end{align*}
\end{thm}

We also get in \cite{nuestro1} a very simple characterization of transient and recurrent Markov chains.

\begin{thm}\label{teo:recurrent_state}
Suppose that both dual Markov chains are irreducible. Then, they are recurrent if and only if the integral
\begin{align*}
\int_\Delta \frac{w_1(x)}{1- {x}} \, \d \mu(x) , && \text{diverges.}
\end{align*}
 Both dual Markov chains are transient whenever the integral converges.
\end{thm}

\begin{cor}[Recurrent and transient Jacobi--Piñeiro Markov chains]\label{pro:recurrentJP}
Both dual Jacobi--Piñeiro Markov chains are recurrent whenever $-1<\alpha_0<0$ and transient if $\alpha_0\geq 0$.
\end{cor}

Applying this result we can see that the Example~\ref{example:JP},
with parameters $\alpha_1 = -\frac 1 4$, $\alpha_2 = \alpha_0 = -\frac 1 2$ is a recurrent Markov chain.
Let see another example of Markov chain associated with the Jacobi--Piñeiros multiple orthogonal polynomials, now with parameters, $ \alpha_1=- \frac{1}{4} $, $\alpha_2 =- \frac{1}{2}$, $\alpha_0 = \frac{1}{2}$, that gives a transient Markov chain.
The corresponding transition matrices,~\eqref{eq:estocastica2} and~\eqref{eq:estocastica1}, are
\begin{align*}
\hat \OH
 \approx
\resizebox{.375\textwidth}{!}{$\left[ \scriptsize\begin{NiceMatrix}
0,3333 & 0.6666& 0 & 0 & 0 &0&0& \Cdots\\ 
0.1026 &0.3205& 0.5769& 0&0 &0&0& \Ddots \\ 
0.0302&0.1163& 0.4712 & 0.3824&0&0&0&\Ddots\\
0& 0.0062 &0.1707& 0.4150 & 0.4080 &0&0& \Ddots\\
0&0&0.0331& 0.1621&0.4600&0.3448&0&\Ddots\\
0&0&0&0.0156&0.1905&0.4279&0.3660 &\Ddots \\
\Vdots&\Ddots&\Ddots&\Ddots&\Ddots&\Ddots&\Ddots&\Ddots
\end{NiceMatrix} \right] $} , 
 &&
\check \OH \approx
\resizebox{.425\textwidth}{!}{$\left[ \begin{NiceMatrix}[columns-width = 0.5cm]
0.3333& 0.3198& 0.3469& 0 & 0 &0&0&0& \Cdots\\ 
0.2138&0.3205& 0.4289& 0.0368& 0 &0&0&0& \Ddots \\ 
0& 0.1565& 0.4711& 0.2726&0.0998&0&0&0&\Ddots\\
0& 0 &0.2395&0.4150& 0.3061&0.0394&0&0& \Ddots\\
0&0&0& 0.2160&0,4600&0.2542&0.0697&0&\Ddots\\
0&0&0&0&0.2583& 0.4279& 0.2746&0.0391&\Ddots \\
\Vdots&\Ddots&\Ddots&\Ddots&\Ddots&\Ddots&\Ddots&\Ddots&\Ddots
\end{NiceMatrix}\right] $} .
\end{align*}
Inspection of both transition matrices, we see that the probabilities to go to the left are bigger in the recurrent situation, with $\gamma=-\frac{1}{2}$, than in the transient case 
with~$\gamma=\frac{1}{2}$.

As we have seen, associated with a tetra diagonal Hessenberg matrix (coming from an AT system) we get two stochastic matrices $\hat \OH$ and $\check \OH$, that we call dual stochastic matrices.
In this sense the corresponding Markov chains are said to be dual.

\begin{thm}
Both dual stochastic matrices, $ \hat \OH$ and $ \check \OH$ are connected by
\begin{align*}
\sigma_{I}^{-1} \, \check \OH \, \sigma_{I} &= \Big( \sigma_{II} \, \hat \OH \, \sigma_{II}^{-1} \Big)^\top
\end{align*}
That is, the stochastic matrices coefficients fulfill
\begin{align*}
\check \OH_{n,n-k}= \frac{B_{n-k} (1) \, Q_{n-k} (1)}{B_{n} (1) \, Q_{n} (1)} \, \hat \OH_{n-k,n},
\end{align*}
with
$m=n-k $, $ k \in \big\{ -1,0,1, 2 \big\} $.
\end{thm}

Now, we will analyze the large $n$ limit for the dual Jacobi--Piñeiro Markov chains.

 \begin{thm} \label{teo:jacobi_pineiro_limit_estocastico}
The large $n$ limit of the Jacobi--Piñeiro stochastic matrices of type~I and II are the same after transposition, i.e.,
\begin{align*}
\lim_{n\to\infty} \check \OH_{n,n+k}
&=\lim_{n\to\infty} \hat \OH_{n+k,n}, & k\in\{-2,-1,0,1\}.
\end{align*}
 \end{thm}
 
\begin{proof}
As we have just seen, we need to show that
\begin{align}\label{eq:other_relations}
\lim_{n \to \infty} \frac{B^{(n-k)}(1) \, Q^{(n-k)}(1)}{B^{(n)}(1) \, Q^{(n)}(1)}= 1, && k=2,1,-1.
\end{align}
%
From~\eqref{eq:Ben1} we directly deduce that
\begin{align} \label{eq:limib}
\lim_{n \to \infty}
\frac{B^{(2n+1)}(1)}{B^{(2n)}(1)} =
\frac{8}{27} =
\lim_{n \to \infty}
\frac{B^{(2n+2)}(1)}{B^{(2n+1)}(1)} .
\end{align}
Now, from Christoffel--Darboux formula for multiple orthogonal polynomials (cf.~\cite{nuestro1})
we get
\begin{multline*}
Q^{(n-1)}(1) \, B^{(n)}(1) = Q^{(n)}(1) \, \big(a_n \, B^{(n-2)}(1) + b_n \, B^{(n-1)}(1) \big)
 \\ +Q^{(n+1)}(1) \, a_{n+1} \, B^{(n-1)}(1),
\end{multline*}
so that we have for the linear forms of type~I the following lower degree homogeneous linear recurrence
\begin{align*}
-Q^{(n-1)}(1) + s_n \, Q^{(n)}(1)+t_n \, Q^{(n+1)}(1)=0,
\end{align*}
with
\begin{align*}
s_n &= a_n \, \frac{B^{(n-2)}(1)}{B^{(n)}(1)}+
b_n \, \frac{B^{(n-1)}(1)}{B^{(n)}(1)}\xrightarrow[n\to\infty]{}
\frac{4^3}{27^3}\frac{27^2}{8^2}+3\frac{4^2}{27^2}\frac{27}{8}=\frac{7}{27},\\
t_n &= a_{n+1} \, \frac{B^{(n-1)}(1)}{B^{(n)}(1)}\xrightarrow[n\to\infty]{}\frac{4^3}{27^3}\frac{27}{8}=\frac{8}{729},
\end{align*}
where we have used~\eqref{eq:limib}, and~\eqref{eq:limit}.
The characteristic polynomial~is
\begin{align*}
-1+\frac{7}{27}r+\frac{8}{729}r^2=\frac{8}{729}(r+27)\Big(r-\frac{27}{8}\Big).
\end{align*}
Therefore, from Poincaré's theorem, having its characteristic roots $\{-27,\frac{27}{8}\}$ distinct absolute value, as the linear forms of type~I are positive at $1$, we get
$\displaystyle \lim_{n\to\infty}\frac{Q^{(n+1)}(1)}{Q^{(n)}(1)}=\frac{27}{8}$
and~\eqref{eq:other_relations} is satisfied.
\end{proof}


We end this work with an example that comes from the hypergeometric weights recently given in~\cite{Lima-Loureiro}.

\begin{exam}
We consider $(W_1, W_2 )$, where the weights are given by
\begin{align*}
 W_1(x) = \omega (x,a,b;c,d), && W_2(x) = \omega (x,a,b+1;c+1,d),
\end{align*}
where $\omega$ is defined on $[0,1]$ as
\begin{align*}
\omega (x,a,b;c,d) =\frac{\Gamma(c)\Gamma(d)}{\Gamma(a)\Gamma(b)\Gamma(\delta)}x^{a-1}(1-x)^{\delta-1}\;
{}_2 {F}{_1}\hspace*{-3pt}\left[{\begin{NiceArray}{c}[small]c-b,d-b \\\delta\end{NiceArray}};1-x\right] ,
\end{align*}
where
$\displaystyle 
{}_2 {F}{_{1}} \left[
\begin{NiceMatrix}[small]a_1, a_2 \\ b_1 \end{NiceMatrix}
;x	\right]
=
\sum_{k=0}^\infty \frac{(a_1)_k \, (a_2)_k}{(b_1)_k }\frac{x^k}{k!}
$
with
$\delta =c+d-a-b$, and
$ a,b,\delta >0$,
$ d-a,d-b,c+1-a,c-b \not\in -\NN \cup \{ 0 \} $.
\end{exam}

Following \cite{Lima-Loureiro} we see that for the hypergeometric tuple
$\big(\frac{4}{3},\frac{5}{3},2,\frac{5}{2}\big)$
we have that
$ c_n=3\kappa$, $b_{n+1}=3\kappa^2$, $a_{n+2}=\kappa^3$, $n \in \NN \cup \{ 0 \}$, with $\kappa = \frac 4{27}$.
That is, the Hessenberg matrix is uniform along its diagonals, i.e. a Toeplitz~matrix
\begin{align}\label{eq:Jacobi_Jacobi_Pineiro_uniform}
\OH =
\resizebox{.3\textwidth}{!}{$\displaystyle
\left[\begin{NiceMatrix}
 3\kappa& 1 & 0 & \Cdots & & \\[-5pt]
 3\kappa^2&3\kappa& 1 & \Ddots& & \\ 
 \kappa^3& 3\kappa^2& 3\kappa& 1 &&\\ 
 0& \kappa^3&3\kappa^2& 3\kappa& 1 & \\
 \Vdots&\Ddots&\Ddots&\Ddots&\Ddots&\Ddots
 \end{NiceMatrix}\right]$} .
\end{align}
We refer to this matrix as the \emph{asymptotic uniform Hessenberg matrix}.
The AT-system $(W_1,W_2,\d x)$ is for this uniform case
\begin{align*} 
W_1(x) & = \frac{81\sqrt 3}{16\pi}\sqrt[3]{x}\left(\sqrt[3]{1+\sqrt{1-x}}-\sqrt[3]{1-\sqrt{1-x}}\right),
 \\
W_2(x) & = \frac{243\sqrt 3}{160\pi}\sqrt[3]{x}\left(\left(\sqrt[3]{1+\sqrt{1-x}}\,\right)^4-\left(\sqrt[3]{1-\sqrt{1-x}}\,\right)^4\right) .
\end{align*}
Note that in~\cite{nuestro2} we find twelve uniform cases.
We also prove in~\cite{nuestro2} a Theorem like~\ref{teo:jacobi_pineiro_limit_estocastico} for the hypergeometric scenario.
Moreover, we could prove that both dual hypergeometric Markov chains are recurrent whenever $0<\delta\leq1$ and transient for $\delta> 1$. In this uniform case $\delta 
 = \frac 7 2$, hence is transient.

We could prove that for some of the twelve uniform tuples (cf.~\cite{nuestro2})
 \begin{align*}
 Q_{n}(1)=\frac 1 {(2\kappa)^n} , && n \in \NN \cup \{ 0 \} .
 \end{align*}
Using this normalization, together with $B_{n}(1) = (2\kappa)^n $, $n \in \NN \cup \{ 0 \} $, we succeed in finding the following semi-stochastic dual matrices associated with $\OH$ in~\eqref{eq:Jacobi_Jacobi_Pineiro_uniform}, and using~\eqref{eq:estocastica2} and~\eqref{eq:estocastica1}, are
\begin{align*}
 \hat \OH = 
 \resizebox{.35\textwidth}{!}{$\displaystyle
 \left[\begin{NiceMatrix}[columns-width =auto]
 \frac{12}{27} &\frac{6}{27} &\frac{1}{27}& 0 & \Cdots& & \\
 \frac{8}{27} & \frac{12}{27} &\frac{6}{27} &\frac{1}{27} & \Ddots & & \\
 0& \frac{8}{27} & \frac{12}{27} &\frac{6}{27} &\frac{1}{27} & \Ddots & \\
 \Vdots& \Ddots& \frac{8}{27} & \frac{12}{27} &\frac{6}{27} &\frac{1}{27} & \\
 & & & \Ddots& \Ddots & \Ddots&\Ddots
 \end{NiceMatrix}
 \right] $},
 && 
 \check \OH
 =
 \resizebox{.285\textwidth}{!}{$\displaystyle
 \left[\begin{NiceMatrix}
 \frac{12}{27}& \frac{8}{27} & 0 & \Cdots& & \\ 
 \frac {6} {27}& \frac{12}{27}& \frac{8}{27} & \Ddots& & \\[5pt]
 \frac{1}{27} & \frac {6} {27}& \frac{12}{27}& \frac{8}{27} &&\\ 
 0& \frac{1}{27} &\frac {6} {27}& \frac{12}{27}& \frac{8}{27} & \Ddots\\ 
 \Vdots&\Ddots&\Ddots&\Ddots&\Ddots&\Ddots
 \end{NiceMatrix}\right]
 $} ,
\end{align*}
which are again Toeplitz Hessenberg matrices.

\section*{Acknowledgments}

The first author thanks Centre for Mathematics of the University of Coim\-bra--UIDB/00324/2020 (funded by the Portuguese Government through FCT/MCTES). The second and third authors acknowledges Center for Research \& Development in Mathematics and Applications 
is supported \linebreak
through the Portuguese Foundation for Science and Technology
(FCT--Fun\-da\-ção para a Ciência e a Tecnologia), references UIDB/04106/2020 and UIDP/04106/2020. The fourth author was partially supported by the Spanish ``Agencia Estatal de Investigación'' research project [PGC2018-096504-B-C33], \emph{Ortogonalidad y~Apro\-ximación: Teoría y Aplicaciones en Física Mate\-má\-tica} and research project [PID2021- 122154NB-I00], \emph{Ortogonalidad y aproximación con aplicaciones en machine learning y teoría de la probabilidad}.


\begin{thebibliography}{99}


\bibitem{afm}
Carlos Álvarez-Fernández, Ulises Fidalgo, and Manuel Mañas,
\emph{Multiple orthogonal polynomials of mixed type: Gauss--Borel factorization and the multi-component 2D Toda hierarchy},
Advances in Mathematics~\textbf{227} (2011) 1451--1525.




\bibitem{abv}
Alexander I. Aptekarev, Amílcar Branquinho, and Walter Van Assche,
\emph{Multiple orthogonal polynomials for classical weights},
Transactions of the American Mathematical Society \textbf{355} (2003) 3887--3914.




\bibitem{Aptekarev_Kaliaguine_VanIseghem}
Alexander Aptekarev, Valery Kaliaguine, and Jeannette Van Iseghem,
\emph{The Genetic Sums's Representation for the Moments of a System of Stieltjes Functions and its Application}, Constructive Approximation \textbf{16} (2000) 487--524.










\bibitem{nuestro2}
Amílcar Branquinho, Juan E. Fernández-Díaz, Ana Foulquié-Moreno, and Manuel Mañas, 
\emph{Hypergeometric Multiple Orthogonal Polynomials and Random Walks},
\hyperref{https://arxiv.org/pdf/2107.00770.pdf}{}{}{\texttt{arXiv:2107.00770}}.





\bibitem{espectral}
Amí­lcar Branquinho, Ana Foulquié-Moreno, Manuel Mañas,
\emph{Spectral theory for bounded banded matrices with positive bidiagonal factorization and mixed multiple orthogonal polynomials},
\hyperref{https://arxiv.org/pdf/2212.10235.pdf}{}{}{\texttt{arXiv:2212.10235}}.




\bibitem{nuestro1}
Amí­lcar Branquinho, Ana Foulquié-Moreno, Manuel Mañas, Carlos Álvarez-Fernández, Juan E. Fernández-Dí­az,
\emph{Multiple Orthogonal Polynomials and Random Walks},
\hyperref{https://arxiv.org/pdf/2103.13715.pdf}{}{}{\texttt{arXiv:2103.13715}}.






\bibitem{Chihara}
Theodore S. Chihara,
\emph{An Introduction to Orthogonal Polynomials},
Gordon \& Breach, 1978, New York. Reprinted by Dover, 2011.




\bibitem{Deift}
Percy Deift, 
\emph{Orthogonal polynomials and random matrices: a Riemann-Hilbert approach},
American Mathematical Soc. Vol. 3, 1999.




\bibitem{Fallat-Johnson}
Shaun M. Fallat and Charles R. Johnson,
\emph{Totally Nonnegative Matrices}, Princeton Series in Applied Mathematics, Princeton University Press, 2011, Princeton.


\bibitem{Gallager}
Robert G. Gallager,
\emph{ Stochastic Processes. Theory for Application}s, Cambridge University Press, 2013.



\bibitem{Gantmacher-Krein}
Felix P. Gantmacher and Mark G. Krein,
\emph{Oscillation and Kernels and Small Vibrations of Mechanical Systems},
revised second edition, AMS Chelsea Publishing, American Mathematical Society, Providence, Rhode Island.



\bibitem{Grunbaum_Iglesia}
F. Alberto Gr\"unbaum and Manuel D. de la Iglesia,
\emph{An urn model for the Jacobi-Piñeiro polynomials},
Proceedings of the American Mathematical Society, 2022,
\hyperref{https://doi.org/10.1090/proc/15910}{}{}{DOI:10.1090/proc/15910}.



\bibitem{Ismail}
Mourad E. H. Ismail,
\emph{Classical and Quantum Orthogonal Polynomails in One Variable},
Encyclopedia of Mathematics and its Applications \textbf{98}, Cambridge University Press,~2009.


\bibitem{KmcG}
Samuel Karlin and James McGregor,
\emph{The Classification of Birth and Death Processes},
Transactions of the American Mathematical Society, \textbf{86} (2) 366--400.




\bibitem{Karlin-McGregor}
Samuel Karlin and James McGregor,
\emph{Random walks},
llinois Journal of Mathematics \textbf{3} (1959) 66--81.



\bibitem{Lima-Loureiro}
Hélder Lima and Ana Loureiro,
\emph{Multiple orthogonal polynomials with respect to Gauss' hypergeometric function},
Studies in Applied Mathematics 2021,
\hyperref{https://doi.org/10.1111/sapm.1243}{}{}{DOI:/10.1111/sapm.12437}.



\bibitem{nikishin_sorokin}
Evgenii M. Nikishin and Vladimir N. Sorokin,
\emph{Rational Approximations and Orthogonality},
Translations of Mathematical Monographs, \textbf{92},
American Mathematical Society, Providence, 1991.


\bibitem{Simon}
Barry Simon,
\emph{Operator Theory. A Comprehensive Course in Analysis, Part 4},
American Mathematical Society, Providence, Rhode Island, 2015.



\bibitem{VanAssche2}
Walter Van Assche,
\emph{Multiple Orthogonal Polynomials}
in Mourad E. H. Ismail,
\emph{Classical and Quantum Orthogonal Polynomials in one Variable},
Encyclopedia of Mathematics and its Applications,
\textbf{98} Cambridge University Press, 2005.
















\end{thebibliography}
\end{document}